\documentclass[12pt]{article}

\usepackage{amsfonts}
\usepackage{amsmath}
\usepackage{amssymb}
\usepackage{indentfirst}

\usepackage[russian]{babel}

\renewcommand{\le}{\leqslant}
\renewcommand{\ge}{\geqslant}

\begin{document}
\title{О хроматических числах некоторых дистанционных графов}
\author{Д. Захаров}
\date{}
\maketitle
\newtheorem{theorem}{Теорема}
\renewcommand{\refname}{Список литературы}

\section{Введение и формулировки результатов}

Положим $\mathcal{R}_n = \{0, 1, \ldots, n-1\}$ и рассмотрим граф
$$
G(n, r, s)=(V, E), \, V=\{v\subset {\cal R}_n: \, |v|=r\}, \, E=\{(v, u): |v\cap u|=s\}.
$$
Этот граф играет большую роль в задачах теории кодирования (см. \cite{WS}), теории Рамсея (см. \cite{Nagy}, \cite{FW}) 
и комбинаторной геометрии (см. \cite{Rai1}--\cite{Rai7}). Именно в связи с последними задачами этот граф интерпретируется как дистанционный.
Особенно важно {\it хроматическое число} этого графа, т.е. величина 
$ \chi(G(n,r,s)) $, равная наименьшему числу цветов, в которые можно так покрасить все вершины графа, чтобы концы любого ребра имели разные цвета. В статьях 
\cite{BKK}, \cite{P} можно найти подробные обзоры всех известных к настоящему времени оценок. В частности, ранее было доказано, что 
$$
n \ge \chi(G(n, r, r-1)) \ge n-r+1
\eqno{(1)}
$$
и что при $ r < 2s+1 $ выполнено неравенство
$$
\chi(G(n, r, s)) \le (1+o(1)) n^{r-s} \frac{r!}{s! ((r-s)!)^2}.
\eqno{(2)}
$$

Нам удалось существенно проще доказать в общем случае оценку (1), нежели это было сделано в работе \cite{BKK}. Также в значительном числе случаев нам удалось 
улучшить эту оценку. А в случае, когда $ r=3, s=2 $, мы нашли бесконечную последовательность значений $ n $, для которых хроматическое число считается точно.
Это первый нетривиальный подобный случай, кроме случая графов $ G(n,3,1) $, для которых точное значение хроматического числа найдено при всех $ n = 2^k $ (см. 
\cite{Rai8}). Наконец, мы смогли улучшить оценку (2). Ниже мы приводим соответствующие формулировки.

\begin{theorem}\label{t1}
Пусть $p>3$ --- простое число. Если для каждого $r$ не выполнено сравнение $-1 \equiv 2^r \pmod p$, 
то для $n = p + 2$ верно неравенство $\chi(G(n, 3, 2))\le n - 2$, а для $n = p + 1$ справедлива оценка $\chi(G(n, 3, 2))\le n - 1$.
\end{theorem}

Теорему \ref{t1} мы докажем в разделе 2, а в разделе 3 мы поясним, почему условиям теоремы удовлетворяет бесконечно много простых чисел $ p $.

\begin{theorem}\label{t2} Выполнены следующие утверждения.

А) Для любых $n, r$ верны неравенства
$$
n\ge \chi(G(n, r, r-1))\ge n-r+1,
$$
если $n-r$ четно, и 
$$
n\ge \chi(G(n, r, r-1))\ge n-r+2,
$$
если $ n-r $ нечетно.

Б) При $n=rk-1$ и простом $r$ верна оценка $\chi(G(n, r, r-1))\ge n-r+2$.
\end{theorem}

В случае $r = 3$ получаем из теоремы \ref{t2}, что
$$
\chi(G(n, 3, 2))\ge n-1, ~~ n\equiv 0, 2, 4, 5 \pmod 6.
$$

Таким образом, для бесконечно многих $n$ найдено значение $\chi(G(n, 3, 2))$. Теорему \ref{t2} мы докажем в разделе 4, а в разделе 5 мы докажем 
следующий результат.

\begin{theorem}\label{t3} Пусть $ n $ --- простое число. Тогда
$$
\chi(G(n, r, s)) \le n^{r-s}.
$$
\end{theorem}

Эту теорему мы докажем двумя способами: с использованием аддитивной комбинаторики и без нее. 

Отметим, что теорема \ref{t3} верна при любых $ r, s $. В то же время из теоремы \ref{t3} нетрудно вывести асимптотическое утверждение, не использующее простоту $ n $,
но слегка ограничивающее параметры $ r, s $: 
$$
\chi(G(n, r, s)) \le (1+o(1))n^{r-s},
$$
коль скоро $ r $ и $ s $ таковы, что при $ n \to \infty $ верна асимптотика $ \left(1+\frac{1}{n^{0.475}}\right)^{r-s} \sim 1 $. Дело в том, что в работе 
\cite{BHP} доказано существование такой константы $ c $, с которой между $ x $ и $ x + cx^{0.525} $ всегда есть простое число. Пользуясь этим результатом, 
получаем
$$
\chi(G(n, r, s)) \le \chi(G(n', r, s)) \le \left(n+cn^{0.525}\right)^{r-s} = 
n^{r-s} \left(1+\frac{c}{n^{0.475}}\right)^{r-s} = (1+o(1))n^{r-s}.
$$
С помощью гипотезы Римана показатель 0.525 умеют заменять на $ 0.5 + o(1) $, и есть вера в то, что, на самом деле, даже на отрезках вида 
$ \left[x,x+c\ln^2 x\right] $ всегда есть простые числа. Разумеется, все это потенциально уменьшает ограничения на $ r $ и $ s $.

Теорема \ref{t3} улучшает прежний результат (2), например, в случае, когда $d=r-s<\sqrt{s}$ и $ s $  достаточно велико. В самом деле, 
используя неравенство $x! \le ex\left (\frac{x}{e}\right )^x$, получаем в (2)
$$
\frac{r!}{s!(d!)^2}=\frac{r(r-1)\ldots(s+1)}{(d!)^2} > \frac{s^d}{(d!)^2} \ge \frac{e^{2d-2}s^d}{d^{2d+2}}. 
$$
Если $ d > \ln s $, то последняя велична по порядку больше, чем 
$$
\frac{s^{d+2}}{s^{d+1}} = s \gg 1.
$$
Иначе та же величина заведомо превосходит 
$$
\left(\frac{s}{\ln^2 s}\right)^d \gg 1.
$$

\section{Доказательство теоремы ~\ref{t1}}

\subsection{Построение раскраски}

\newtheorem{lemma}{Лемма}

Сразу заметим, что достаточно рассмотреть случай $n = p + 2$. В самом деле, чтобы построить требуемую раскраску для графа 
$G(n, 3, 2)$ в $n-1$ цвет, достаточно взять раскраску графа $G(n+1, 3, 2)$ в $n-1$ цвет и убрать все вершины (тройки элементов), проходящие через какой-то элемент. 

Итак, построим некоторую раскраску вершин нашего графа. Выделим элементы $n-1$ и $n-2$ и разобьем множество вершин графа $G(n, 3, 2)$ на 4 подмножества:

$V_0$ --- множество вершин, непересекающихся с $n-1$ и $n-2$,

$V_2$ --- множество вершин, пересекающихся с $n-1$ и $n-2$,

$W_1$ --- множество вершин, непересекающихся с $n-1$ и пересекающихся с $n-2$,

$W_2$ --- множество вершин, пересекающихся с $n-1$ и непересекающихся с $n-2$.

\noindent Каждую вершину  $x=(x_1, x_2, x_3)\in V_0$ покрасим в цвет 
$$
\chi(x)=x_1+x_2+x_3 \pmod p,
$$
 каждую вершину $x=(x_1, n-1, n-2)\in V_2$ покрасим в цвет 
 $$
 \chi(x)=3x_1 \pmod p.
$$
Для раскраски множеств $W_1$ и $W_2$ воспользуемся следующей леммой, которую мы докажем в параграфе 2.3.

\begin{lemma}\label{l}

Если обозначить $ \binom{\mathcal{R}_{p}}{2} $ множество всех пар элементов из $ {\cal R}_p = {\cal R}_{n-2} $, то существуют такие
функции $f_1$, $f_2$: $\binom{\mathcal{R}_{p}}{2} \rightarrow \mathcal{R}_{p}$, что

\begin{enumerate}

\item  $f_i(x, y)=f_i(y, x) \in \{x, y\}$;

\item $f_1(x, y) \not = f_2(x, y)$;

\item не существует такой пары $ x, y $, для которой при каком-нибудь $i \in \{1, 2\}$ было бы выполнено
$$
f_i(x, y)=x, ~~ f_i\left(\frac{x+y}{2}, x\right) = \frac{x+y}{2} 
$$
(здесь деление --- это деление в $\mathbb{Z}_p$, т.е. выбор соответствующего вычета в множестве $ \{0, \dots, p-1\} $).

\end{enumerate}

\end{lemma}

Теперь покрасим каждую вершину 
$$
x=(x_1, x_2, m)\in W_i, ~~ m\in\{n-1, n-2\},
$$
в цвет
$$
\chi(x)=x_1+x_2+f_i(x_1, x_2) \pmod p.
$$ 

Очевидно, что мы построили раскраску всех вершин нашего графа в $ p = n-2 $ цвета. В следующем параграфе мы докажем, что полученная раскраска 
является правильной, т.е. что любые две вершины (тройки элементов), пересекающиеся ровно по двум элементам (образующие ребро), покрашены в разные цвета.  

\subsection{Доказательство правильности раскраски}

Возьмем любые две смежные вершины 
$$
x = (x_1, x_2, x_3), ~~~ y = (x_1, x_2, x_4) 
$$
и рассмотрим следующие случаи.

\paragraph{Случай 1.} Пусть $x, y \in V_0$. Тогда, очевидно, 
$$
\chi(x) \equiv x_1+x_2+x_3 \not \equiv x_1+x_2+x_4 \equiv \chi(y) \pmod p.
$$

\paragraph{Случай 2.} Пусть $x, y \in V_2$. Тогда 
$$
\chi(x) \equiv 3x_3\not \equiv 3x_4 \equiv \chi(y) \pmod p,
$$
так как по условию $ p $ простое, большее трех, т.е. $p \not \vdots 3$.

\paragraph{Случай 3.} Пусть $x \in V_0, $ $ y \in V_2$. Этот случай, очевидно, невозможен.

\paragraph{Случай 4.} Пусть $x \in W_i, $ $ y \in V_2$. Тогда можно считать, что $x_1 = n-1 $, $ x_4 = n-2$, откуда
$$
\chi(x)-\chi(y)\equiv x_2+x_3+f_i(x_2, x_3)-3x_2 \equiv x_3 + f_i(x_2, x_3) - 2x_2 \pmod p.
$$ 
Если $f_i(x_2, x_3)=x_2$, то
$$
\chi(x)-\chi(y) \equiv x_3 - x_2 \not \equiv 0 \pmod p,
$$ 
а если $f_i(x_2, x_3)=x_3$, то 
$$\chi(x)-\chi(y) \equiv 2x_3 - 2x_2 \not \equiv 0 \pmod p,$$
так как $p=n - 2$ нечетно.

\paragraph{Случай 5.} Пусть $x\in W_i $, $y \in V_0$. Можно считать, что $x_3 = n-1$. Имеем 
$$
\chi(x)-\chi(y) \equiv x_1+x_2+f_i(x_1, x_2)-x_1-x_2-x_4 \equiv f_i(x_1, x_2)-x_4 \not \equiv 0 \pmod p
$$
по свойству 1 из леммы \ref{l}.

\paragraph{Случай 6.} Пусть $x \in W_1$, $y \in W_2$. В этом случае $x_3=n-2$, а $x_4=n-1$, и, значит,
$$
\chi(x)-\chi(y)\equiv f_1(x_1, x_2)-f_2(x_1, x_2) \not \equiv 0 \pmod p
$$
по свойству 2 леммы \ref{l}. 

\paragraph{Случай 7.} Пусть $x, y \in W_i$. Можно считать, что $i=1$, то есть $x_1 = n-2$. Можем написать
$$
\chi(x)-\chi(y) \equiv x_2+x_3+f_i(x_2, x_3)-x_2-x_4-f_i(x_2, x_4)\equiv x_3+f_i(x_2, x_3)-x_4-f_i(x_2, x_4) \pmod p. 
$$
Рассмотрим подслучаи.

\paragraph{Подслучай 7.1.} Если $f_i(x_2, x_3)=f_i(x_2, x_4)=x_2$, то
$$
x_3+f_i(x_2, x_3)-x_4-f_i(x_2, x_4)\equiv x_3-x_4 \not \equiv 0 \pmod p.
$$

\paragraph{Подслучай 7.2.} Если $f_i(x_2, x_3)=x_3, $ $ f_i(x_2, x_4)=x_4$, то
$$
x_3+f_i(x_2, x_3)-x_4-f_i(x_2, x_4)\equiv 2x_3-2x_4 \not \equiv 0 \pmod p.
$$

\paragraph{Подслучай 7.3.} Если $f_i(x_2, x_3)=x_2, $ $ f_i(x_2, x_4)=x_4$, то
$$
x_3+f_i(x_2, x_3)-x_4-f_i(x_2, x_4)\equiv x_3+x_2-2x_4 \pmod p.
$$
Предположим, что $x_3+x_2\equiv 2x_4 \pmod p$. Имеем $x_4=\frac{x_2+x_3}{2}$. Но тогда мы получаем, что 
$$
f_i(x_2, x_3)=x_2, ~~~ f_i\left(x_2, \frac{x_2+x_3}{2}\right)=\frac{x_2+x_3}{2}, 
$$
что противоречит свойству 3 леммы \ref{l}. Значит, $\chi(x)-\chi(y)\not\equiv 0 \pmod p$.

\vskip+0.3cm

Таким образом, мы разобрали все случаи, и раскраска является правильной.

\subsection{Доказательство леммы \ref{l}}

Для каждой пары различных $i$ и $j$ построим следующие циклические последовательности (далее будем называть их окружностями) в $ \mathbb{Z}_p$: 
$$
C(i, j)=\left(j_m \left|j_0=j, ~~~ j_m\equiv \frac{j_{m-1}+i}{2} \right. \pmod p\right).
$$
Например, для $p = 7$ такой окружностью будет 
$$
C(1, 5)=C(1,3)=C(1, 2)=(2, 5, 3)=(3,2,5)=(5,3,2).
$$
\begin{lemma}
Длина $|C(i, j)|$ окружности $C(i, j)$ равна порядку $k$ числа $2$ в $\mathbb{Z}_p$. 
\end{lemma}

\paragraph{Доказательство.} Выразим $j_m$ через $i$ и $j$:
$$
j_m\equiv\frac{j_{m-1}+i}{2}\equiv\frac{\frac{j_{m-2}+i}{2}+i}{2}\equiv\frac{j_{m-2}+3i}{4}\equiv\frac{\frac{j_{m-3}+i}{2}+3i}{4}\equiv\frac{j_{m-3}+7i}{8}
$$
Ясно, что 
$$
j_m\equiv\frac{j+(2^m-1)i}{2^m}  \pmod p.
\eqno{(3)}
$$
Но тогда
$$
j_k\equiv\frac{j+(2^k-1)i}{2^k}\equiv i+\frac{j-i}{2^k}\equiv i+j-i \equiv j_0 \pmod p.
$$
Если же предположить, что при $m<k$ выполнено $j_m \equiv j_0 \pmod p$, то
$$
2^m(j_0-i) \equiv j_0-i \pmod p, 
$$
$$
(2^m - 1)(j_0-i)\equiv 0 \pmod p,
$$
что неверно для простых $p$. Лемма доказана.

\vskip+0.3cm

Построим граф $G=(V, E)$ окружностей следующим образом: вершины --- все окружности $C(i, j)$, а ребра --- пары окружностей 
$C(i, j)$ и $C(j, i)$. Мы соединяем именно циклические последовательности: 
например, для $p=7$ окружности $C(1, 2), C(1, 3), C(1, 5)$ представляют одну и ту же вершину графа, которая,
тем самым, соединена с каждой из вершин 
$$ 
C(2,1) = C(2,5) = \ldots, ~~~ C(3,1) = \dots, ~~~ C(5,1) = \ldots, 
$$
и т.д.

\begin{lemma} 
$\chi(G)=2$.
\end{lemma}

\paragraph{Доказательство.} Предположим, что в $G$ нашелся цикл нечетной длины $l$:
$$
C(i_1, i_2), ~ C(i_2, i_3), ~ \ldots, ~ C(i_{l-1}, i_l), ~ C(i_l, i_1).
$$
Ясно, что наличие ребра равносильно тому, что $i_{m-1} \in C(i_{m}, i_{m+1})$. Согласно (3) для некоторого $s_m$ имеем
$$
i_{m-1}\equiv \frac{i_{m+1}+(2^{s_m}-1)i_m}{2^{s_m}} \pmod p,
$$
$$
-2^{s_m}(i_m-i_{m-1})\equiv i_{m+1}-i_m  \pmod p.
$$
Пусть $d_m = i_m-i_{m-1}$. Тогда последнее равенство можно переписать в виде
$$
d_{m+1}\equiv -2^{s_m}d_m \pmod p.
$$
Применяя эту формулу $l$ раз, получаем
\begin{eqnarray*}
&d_l \equiv -2^{s_{l-1}}d_{l-1} \equiv (-1)^2 2^{s_{l-1}+s_{l-2}}d_{l-2} \equiv \ldots \qquad&\\  
&\qquad \ldots \equiv (-1)^{l-1}2^{s_{l-1}+s_{l-2}+\ldots+s_1}d_1\equiv (-1)^l 2^{s_1+\ldots+s_l}d_l. &
\end{eqnarray*}
Значит, для некоторого $S$ выполнено
$$
d_l(2^S+1)\equiv 0 \pmod p,
$$
что невозможно по условию теоремы. Лемма доказана.

\vskip+0.3cm

Перейдем к завершению доказательства леммы \ref{l}, то есть к построению функций, обладающих тремя свойствами из формулировки леммы.
А именно, раскрасим все окружности в два цвета: $\mathcal{M}_1$ --- множество окружностей первого цвета, $\mathcal{M}_2$ --- второго. 
И положим

$f_1(i, j) = i, f_2(i, j) = j$, если $C(i, j)\in \mathcal{M}_1$;

$f_1(i, j) = j, f_2(i, j) = i$, если $C(i, j)\in \mathcal{M}_2$.

Тем самым, свойства 1 и 2 автоматически выполнены. Докажем свойство 3. Пусть нашлись такие $x, y$, что 
$$
f_i(x, y)=x, ~~~ f_i\left(\frac{x+y}{2}, x\right)=\frac{x+y}{2}.
$$
Тогда
$$
C(x, y) \in \mathcal{M}_i,~~~ C\left(\frac{x+y}{2}, x\right)\in \mathcal{M}_i,
$$
а значит,
$$
C\left(x, \frac{x+y}{2}\right)\in \mathcal{M}_{3-i},
$$
но ведь ясно, что $C(x, y)=C\left(x, \frac{x+y}{2}\right)$. Приходим к противоречию, и лемма доказана.

\section{Доказательство бесконечности числа простых \\ в теореме \ref{t1}}

Докажем, что условию $-1\not\equiv 2^r \pmod p$ удовлетноряют все простые числа вида $p=8k-1$, которых бесконечно много по 
классической теореме Дирихле. Заметим, что условию удовлетворяют не только они. Например, 73 тоже удовлетворяет условию.

Пусть $d$ --- порядок двойки по модулю $p=8k-1$. Обозначим $ \left(\frac{a}{p}\right) $ символ Лежандра (см. \cite{Vin}). Хорошо известно, что (см. \cite{Vin}) 
$$
\left(\frac{2}{p}\right)=(-1)^{\frac{p^2-1}{8}}=(-1)^{\frac{64k^2-16k+1-1}{8}}=1.
$$
Значит, 2 --- квадратичный вычет в $\mathbb{Z}_p$. Следовательно, 
$2^{\frac{p-1}{2}}\equiv 1 \pmod p$, а стало быть, 
$d|\frac{p-1}{2}=4k-1$, то есть $d$ нечетно. 

Предположим, что существует $l$ --- наименьшее натуральное число, такое, что $2^l\equiv -1 \pmod p$.  
Тогда $d | 2l$, но $d$ нечетно, откуда следует, что $d|l$ и, значит, $2^l \equiv 1 \pmod p$. Полученное противоречие завершает доказательство.

\section{Доказательство теоремы ~\ref{t2}}
Сначала установим неравенство $\chi(G(n, r, r-1))\le n$. Сопоставим каждой вершине $x=(x_1, \ldots, x_r)$ цвет 
$$
\chi(x)=x_1+\ldots+x_r \pmod n.
$$
Тогда для двух смежных вершин $x=(x_1, \ldots, x_r)$ и $y=(x_2, \ldots, x_{r+1})$ разность цветов будет равна $\chi(x)-\chi(y)\equiv x_1-x_{r+1}\not \equiv 0 \pmod n$.

Перейдем к нижним оценкам из пункта А). 
Пусть $W$ --- независимое множество в $G(n, r, r-1)$. Тогда зафиксируем $r-2$ элемента $a_1, \ldots, a_{r-2}$. 
Понятно, что их содержит не более $[\frac{n-r+2}{2}]$ элементов множества $W$, а каждое множество из $W$ содержит $\binom{r}{r-2}$ множеств мощности $r-2$, откуда
$$
|W|\le \frac{[\frac{n-r+2}{2}]\binom{n}{r-2}}{\binom{r}{2}},
$$
то есть 
$$
\chi(G(n, r, r-1))\ge\frac{|V(n, r, r-1)|}{\alpha(G(n, r, r-1))}\ge
\frac{\binom{n}{r}\binom{r}{2}}{[\frac{n-r+2}{2}]\binom{n}{r-2}}=\frac{(n-r+2)(n-r+1)}{2[\frac{n-r+2}{2}]}.
$$
Рассматривая два случая четности разности, получаем заявленное в пункте А) теоремы.

Пусть, наконец, $n=rk-1$, $ r $ простое. Предположим, что вершины графа 
$G(n, r, r-1)$ покрашены в $n-r+1$ цвет. Тогда в каждой максимальной клике встречаются все цвета. 
Всего максимальных клик $\binom{n}{r-1}$, а каждая вешина содержится в $r$ кликах. Значит, вершин каждого цвета 
$\frac{\binom{n}{r-1}}{r}$, но, так как $r$ простое и $n=rk-1$, то это число нецелое. Противоречие. 

Теорема доказана.

\section{Доказательство теоремы ~\ref{t3}}

\subsection{Первое доказательство}

Пусть $n$ простое, $h=r-s$. Тогда по теореме Бозе--Човла (см. \cite{HR}) 
существуют такие числа $a_0, \ldots, a_{n-1} \in \mathbb{Z}_{n^h-1}$, что все суммы по $h$ чисел из этого набора различны. 

Теперь строим раскраску следующим образом: если $x=(x_1, \ldots, x_r) $ --- вершина, то 
$$
\chi(x)=a_{x_1}+\ldots+a_{x_r} \pmod {n^h-1}.
$$
Понятно, что тогда общее количество цветов не больше $n^h$. Проверим правильность расскраски. Действительно, если вершины
$$
x=(x_1, \ldots, x_s, y_1, \ldots, y_{r-s}), ~~~ y=(x_1, \ldots, x_s, z_1, \ldots, z_{r-s})
$$ 
образуют ребро, то 
$$
\chi(x)-\chi(y) \equiv a_{y_1}+\ldots+a_{y_h}-\left(a_{z_1}+\ldots+a_{z_h}\right) \not \equiv 0 \pmod {n^h-1}
$$
по определению $a_i$.

Теорема доказана.

\vskip+0.3cm

Отметим, что в построенной раскраске каждый цвет состоит из множеств, мощности пересечения которых не только не равны $ s $, но строго меньше $ s $. 
Это позволяет рассчитывать на дальнейшее усиление результата. 

Отметим также, что в теореме Бозе--Човла можно брать не только простые числа, но и степени простых. Однако это мало влияет на общность результата (ср. 
замечания после формулировки теоремы, в которых говорится о плотности распределения простых в натуральном ряде).

\subsection{Второе доказательство}
Пусть снова $n$ простое, $h=r-s$.
Цветами у нас будут векторы в $\mathbb{Z}_n^h$. А именно, пусть $x=(x_1, \ldots, x_r)$. Тогда $i$-ю координату цвета определим как сумму всех 
произведений из $i$ координат $x$ по модулю $n$. Предположим, что раскраска не является правильной, то есть найдутся 
$x=(x_1,\ldots,x_h, z_1, \ldots, z_s)$ и  $y=(y_1,\ldots,y_h, z_1, \ldots, z_s)$ одного цвета. 
Положим 
$$
\mathbb{X}=\{x_1, \ldots, x_h\}, ~~ \mathbb{Y}=\{y_1, \ldots, y_h\}, ~~ \mathbb{W}=\{z_1,\ldots, z_s\},
$$ 
а также 
$$
\sigma_i (\mathbb{M})=\sum_{S\subset\mathbb{M}, |S|=i}\prod_{t\in S}t, ~~~ \mathbb{M} \in \{\mathbb{X}, \mathbb{Y}, \mathbb{W}\}.
$$
Тогда докажем по индукции, что для любого $i \in \{1,\ldots, h\}$
$$
\sigma_i(\mathbb{X})\equiv \sigma_i(\mathbb{Y}) \pmod n.
$$

\paragraph{База.} При $i = 1$ условие одноцветности вершин $ x, y $ означает, что
$$
\sum_{i=1}^h x_i + \sum_{i=1}^s z_i \equiv \sum_{i=1}^h y_i + \sum_{i=1}^s z_i \pmod n,
$$
то есть $\sigma_1(\mathbb{X}) \equiv \sigma_1(\mathbb{Y}) \pmod n$, что и требовалось.
 
\paragraph{Переход.} Пусть для всех $i_0<i$ утверждение верно. Тогда если мы рассмотрим суммы слагаемых, в каждом из которых 
$i$ множителей, то слагаемые, представляющие собой произведения только элементов множества $\mathbb{W}$, сократятся сразу. 
В свою очередь, суммы тех слагаемых, в каждое из которых входит ровно по $j < i$ множителей из $\mathbb{W}$, представятся в виде 
$\sigma_{i-j}(\mathbb{M})\sigma_j(\mathbb{W})$, где
$\mathbb{M} \in \{\mathbb{X}, \mathbb{Y}\}$. Таким образом, по предположению индукции эти суммы также сократятся, а значит, 
сократятся все слагаемые, кроме тех, в которых нет элементов множества $\mathbb{W}$, откуда
$$
\sigma_i(\mathbb{X})\equiv \sigma_i(\mathbb{Y}) \pmod n.
$$
Переход сделан.

Таким образом, мы получаем, что для множеств $\mathbb{X}, \mathbb{Y}$, имеющих мощность $ h $ каждое, выполнены сравнения
$$\sigma_1(\mathbb{X})\equiv \sigma_1(\mathbb{Y}) \pmod n,$$
$$\sigma_2(\mathbb{X})\equiv \sigma_2(\mathbb{Y})\pmod n,$$
$$ \vdots $$
$$\sigma_h(\mathbb{X})\equiv \sigma_h(\mathbb{Y})\pmod n.$$
Но тогда по теореме Виета эти множества являются корнями некоторого многочлена степени $h$, а значит, они совпадают, что невозможно. 

Теорема доказана.

\end{document}